\documentclass[10pt]{amsart} \usepackage{amsfonts}
\usepackage{color}
\RequirePackage[OT1]{fontenc}
\RequirePackage{amsthm,amsmath}
\RequirePackage[numbers,square]{natbib}

\definecolor{c20}{rgb}{0.,0.7,0.}
\definecolor{c30}{rgb}{0.,0.,1.}
\definecolor{c40}{rgb}{1,0.1,0.7} \definecolor{c50}{rgb}{1,0,0}
\definecolor{c60}{rgb}{1,0.9,0.1}

\def\cL#1{\textcolor{c50}{#1}}  
\def\cL#1{#1}

\newcommand{\kb}[1]{\boldsymbol{#1}}
\newcommand{\vk}[1]{\kb{#1}} 

\newcommand{\ve}{\varepsilon}
 
\newcommand{\abs}[1]{\left\lvert #1 \right\rvert}

\newcommand{\E}[1]{\mathbb{E}\left\{#1\right\}}

\newcommand{\pk}[1]{\mathbb{P} \left\{ #1 \right \} }

\newcommand{\R}{\mathbb{R}}

\newcommand{\N}{\mathbb{N}} 
\newcommand{\inn}{\in \N} \newcommand{\ldot}{,\ldots,}

\newcommand{\BQN}{\begin{eqnarray}}
\newcommand{\EQN}{\end{eqnarray}}
\newcommand{\BQNY}{\begin{eqnarray*}}
\newcommand{\EQNY}{\end{eqnarray*}}

\newcommand{\BS}{\begin{sat}} \newcommand{\ES}{\end{sat}}
\newcommand{\BT}{\begin{theo}} \newcommand{\ET}{\end{theo}}
\newcommand{\BL}{\begin{lem}} \newcommand{\EL}{\end{lem}}
\newcommand{\BK}{\begin{korr}} \newcommand{\EK}{\end{korr}}

\newcommand{\BD}{\begin{de}} \newcommand{\ED}{\end{de}}
\newcommand{\BIT}{\begin{itemize}}
\newcommand{\EIT}{\end{itemize}}
\newcommand{\BDI}{\begin{description}}
\newcommand{\EDI}{\end{description}}

\newcommand{\BRM}{\begin{remarks}}
\newcommand{\ERM}{\end{remarks}}

\newcommand{\BEL}{\begin{lem}} \newcommand{\EEL}{\end{lem}}

\newtheorem{theo}{Theorem}[section]
\newtheorem{sat}[theo]{Proposition}
\newtheorem{de}[theo]{Definition} \newtheorem{lem}[theo]{Lemma}

\newtheorem{korr}[theo]{Corollary}
\newtheorem{remark}[theo]{Remark}
\newtheorem{remarks}[theo]{Remarks}

\newcommand{\nelem}[1]{{Lemma \ref{#1}}}
\newcommand{\neprop}[1]{{Proposition \ref{#1}}}
\newcommand{\netheo}[1]{{Theorem \ref{#1}}}

\newcommand{\prooftheo}[1]{ \textsc{\bf Proof of Theorem}
\ref{#1}:} \newcommand{\proofprop}[1]{\textsc{\bf Proof of
Proposition} \ref{#1}:} 
  \newcommand{\COM}[1]{}

\newcommand{\QED}{\hfill $\Box$}

  \def\IF{\infty}

\newcommand{\expon}[1]{\exp\left(#1\right)}

  \date{}

\def\FRE{Fr$\mathrm{\acute{e}}$chet} \def\I#1{\mathbb{I}\left\{#1>u\right\}}
\def\equaldis{\stackrel{d}{=}}
\newcommand{\todis}{\stackrel{d}{\to}}
\newcommand{\nesec}[1]{{Section \ref{#1}}}

\def\pru{p_r(u)}
\def\tilde{\widetilde}

\begin{document}

\title{Extremes of Order Statistics of Self-similar Processes}

\COM{\author{Patrik Albin}
\address{Department of Mathematical Sciences, Chalmers University of
Technology, Sweden}
\email{palbin@chalmers.se}

\author{Enkelejd  Hashorva}
\address{Enkelejd Hashorva, Department of Actuarial Science,
University of Lausanne,\\
UNIL-Dorigny, 1015 Lausanne, Switzerland

\email{Enkelejd.Hashorva@unil.ch}
}

\author{Lanpeng Ji}
\address{Lanpeng Ji, Department of Actuarial Science,
University of Lausanne\\
UNIL-Dorigny, 1015 Lausanne, Switzerland
}
\email{Lanpeng.Ji@unil.ch}

\author{Chengxiu Ling }
\address{Chengxiu Ling, Department of Actuarial Science,
University of Lausanne,\\
UNIL-Dorigny, 1015 Lausanne, Switzerland }
\email{Chengxiu.Ling@unil.ch}
}

\author{Chengxiu Ling }
\address{Chengxiu Ling, Department of Actuarial Science,
University of Lausanne,\\
UNIL-Dorigny, 1015 Lausanne, Switzerland }
\email{Chengxiu.Ling@unil.ch}
\bigskip

\date{\today} \maketitle
\def\pru{p_r(u)}

{\bf Abstract:} Let $\{X_i(t),t\ge0\}, 1\le i\le n$ be independent
copies of a random process $\{X(t), t\ge0\}$. For a given positive constant
$u$, define the set of $r$th conjunctions $C_r(u):=\{t\in[0,1]: X_{r:n}(t)>u\}$ with  
$ X_{r:n}$ the $r$th largest order statistics of $X_i, 1\le i\le n$.  In numerical applications such as brain mapping and digital communication systems, of interest is the approximation of $\pru=\pk{C_r(u)\neq\phi}$. Instead of  stationary processes dealt with in \cite{KHorders},  we consider in this paper $X$ a self-similar $\R$-valued process with $P$-continuous sample paths. By imposing the Albin\rq{}s conditions directly on $X$,  we establish an exact asymptotic expansion of $\pru$ as $u$ tends to infinity. As a by-product we derive the asymptotic tail behaviour of the
mean sojourn time of $X_{r:n}$ over an increasing threshold. Finally, our findings are illustrated for the case that $X$ is
a bi-fractional Brownian motion, a sub-fractional Brownian motion,
and a generalized self-similar skew-Gaussian process.

\COM{For each index $r\le n$ the process $X_{r:n}$ is also
self-similar, and thus the asymptotics of the tail probability
$\pru=\pk{\sup_{t\in [0,1]} X_{r:n}(t)> u}$ as $u\to \IF$ can be
derived by imposing the Albin's conditions. Since by doing so, for
any $r$, a set of conditions will be imposed,
In this paper we show
that it is possible to derive the asymptotics of $\pru$ as $u\to
\IF$ for any $T>0$ by imposing Albin's condition directly on the process $X$.}



\section{Introduction}\label{sec1} Let $\{X(t),t\ge 0\}$ be a
self-similar $\R$-valued random process with index $\kappa>0$ and $P$-continuous sample
paths, i.e., the
finite-dimensional distributions (f.d.d.) of $X(\lambda t)$
coincide with those of $\lambda^\kappa X(t)$ for all $\lambda
>0$. Given $X_1 \ldot X_n, n\inn$ independent copies of
$X$ we define the $r$th order statistics process $X_{r:n}, 1\le r\le n$ as follows.
\BQN \label{def: Orderp} 
X_{n:n}(t) \le \cdots \le X_{1:n}(t), \quad t>0. 
\EQN
Order statistics play an important role in many statistical applications. 
If $X_i(t)$ models the value of a certain  object (say signal or image) $i$
measured at time $t\in [0,1]$, and $u$ a fixed threshold, then the $r$th conjunction occurs in the time set  $C_r(u):=\{t\in[0,1]: X_{r:n}(t)>u\}$. In applications, of interest is to obtain the probability that $C_r(u)$ is not empty, which is given by 
\BQN\label{TailOS} 
\pru: = \pk{\sup_{t\in [0,1]}X_{r:n}(t)> u}.
\EQN
Clearly, $\pru$ is the probability that
at least $r$ objects overshoot the threshold $u$  before time point 1. Most numerical applications, related to the analysis of functional magnetic resonance
imaging (fMRI) data and the surface roughness during all machinery processes, are concerned with the calculation of 
$\pru$; \cL{see, e.g., \cite{AdlerST2012, MR1747100}.} \\
For certain smooth Gaussian random fields
approximations of $\pru$ are discussed for instance in
\cite{MR2775212, ChengXiao13, MR1747100}; results for
non-Gaussian random fields can be found for instance in
\cite{MR2654766}. The recent contributions
\cite{DebickiHJminima, KHorders} derive asymptotic expansions
of $\pru$ considering a stationary (Gaussian) process $X$. It is
well-known that the stationary random field cannot be used to
model phenomena and data sets that exhibit certain
non-stationary characteristics such as long-range dependence
(LRD). Such situations arise naturally in the limit theorems of
random walks and other stochastic processes modeling various phenomena such as 
telecommunications, internet  traffic, image processing and 
mathematical finance. Typically, these processes are related to self-similar
processes. We refer to the monographs \cite{AlbinS2004, BardetT2010, BeranFHK2013, 
EmbrechtsM2002, Taqqu1978} for complete
expositions on theoretical and practical aspects of
self-similar processes.

Since the exact calculation of $\pru$ is not possible in general, in this contribution we derive the approximations of  $\pru, u\to\IF$ for a self-similar generic process $X$ with index $\kappa$. \\
For the formulation of our main result, we need to introduce Albin\rq{}s conditions imposed on $X$; see, e.g.,  \cite{Albin1998, KHorders}. Hereafter, let $\overline F=1-F$ denote the survival function of a given distribution $F$.

{\bf Condition A}: (Gumbel MDA and weak convergence) Suppose that $X(1)$ has a continuous distribution function (df) $G$ with infinite right endpoint, and it belongs to the Gumbel max-domain attraction (MDA) with some positive scaling function $w=w(u)$, denoted by $G\in MDA(\Lambda, w)$, that is
\BQNY \lim_{u\to\IF}
\frac{\overline G(u+x/w)}{\overline G(u)} = e^{-x}, \quad x\in\R. 
\EQNY
Further, we assume that 
there is a non-increasing, positive function $q=q(u) = D_0 u^{-\alpha_0} (1+o(1)), \, u\to\IF$
with some positive $D_0, \alpha_0$,  \COM{such that
\BQN
\label{eqQ} \lim_{u\to\IF} \frac{1}{ q(u)} =Q \in (0,\IF],
\quad\tilde a = \frac{1}{2\sup_{u \in [0, \IF)}q(u)} >0
 \EQN}
 and an $\overline \R$-valued random process
$\{\xi(t), t\ge 0\}$ with 0 the continuity of the df of $\xi(t), t\ge0$ such that
for all $m\in\mathbb N, t_i\ge0, 1\le i\le m$
\BQN\label{Eq: condL} 
\lim_{u\to\IF}\pk{X(1-qt_i)>u, 1\le i\le m \lvert
X(1)>u} = \pk{\xi(t_i) > 0, 1\le i\le m }. 
\EQN
{\bf Condition B}: (Short-lasting-exceedance) For the function
$q$ as in condition A we have 
\BQNY
\lim_{d\to
\IF}\limsup_{u\to \IF} \int_{d\wedge (1/q)}^{1/q}
\pk{X(1-qt)>u  \lvert X(1) > u}\, dt = 0
\EQNY
with $a\wedge b= \min(a, b)$. \\
In order to state the tightness condition C below
we need the following notation. Let $\tilde a= 1/(2\sup_{u\in[0,\IF)}q(u))$ be positive, and  $L(u)$  the mean sojourn time of the process $X$ over the
threshold $u$ which is defined as 
$$
L(u) = \int_0^1 \I{X(t) }\, dt
$$ 
where $\mathbb{I}\{\cdot\}$ is the indicator function. 
Further we define a  sequence $\{t_a^u(k)\}$ for given $u\in\R, a\in(0, \tilde a]$: setting $t_a^u(0) =1$ and
\BQN\label{t-sequence} 
t_a^u(k+1) = t_a^u(k) \big(1- aq((t_a^u(k))^{-\kappa}u) \big),\quad k=0, \ldots, K 
\EQN
with $ K\equiv K(a, u) = \sup\{k\in \mathbb N: (t_a^u(k))^{-\kappa} u < \IF\}$.

{\bf Condition C}: For the functions $w$ and $q$ given as in condition A, we have for each $\sigma >0$
\BQNY \vk{v}(a, \sigma) \equiv \limsup_{u\to\IF}
\frac{\pk{\sup_{t\in[0,1]} X(t) > u+\sigma/w, \max_{0\le k\le K}
X(t_a^u(k)) \le u } } { \E{L(u)/q} + \overline G(u)} \to0,\quad a\downarrow0. 
\EQNY
Condition C is often verified via Propositions
3--5 in \cite{Albin1998}. For instance, in view of Proposition
3(ii) in \cite{Albin1998} sufficient for condition C to
hold is the next condition.

{\bf Condition C$^*$}: Suppose that there exist some positive
constants $\lambda_0, \rho, b, D$ and $d>1$ such that
 \BQNY
\pk{X(1-qt) > u+\frac{\lambda+v}{ w}, X(1) \le u +\frac vw}
\le D t^d\lambda^{-b} \overline G(u) 
\EQNY 
holds for all $u$ large and all $ 0 < t^\rho \le \lambda \le \lambda_0, v\ge0$.
\COM{
 We abbreviate that as  $G\in MDA(\Lambda,w)$.
\\ Note in passing that (see, e.g., \cite{Faletal2010}) the assumption
$G\in MDA(\Lambda,w)$ implies $\lim_{u\to\IF} uw =\IF$ and
further 
\BQN\label{self-neglecting} 
\frac{w(u+x/w)}{w}
\to1, \quad u\to\IF 
\EQN 
holds locally uniformly for $x\in\R$.\\}

The contribution of this paper is to establish the approximation of $\pru$ under the common  Albin\rq{}s conditions A, B and C (or C$^*$) given above, which are satisfied by many well-known self-similar processes such as skewed $\alpha$-stable processes, Kesten-Spitzer processes and Rosenblatt processes in \cite{Albin1998}. We will illustrate our findings by the following self-similar processes in \nesec{sec4}:
i) a bi-fractional Brownian motion, ii) a sub-fractional
Brownian motion, and iii) a generalized self-similar
skew-Gaussian process, see Theorem \ref{T5}. \cL{We refer to 
\cite{BojdeckiGT2007, BojdeckiGT2012, HuslerP1999, LiX2011, Sato1991, Sghir2014} for related studies on these generalized self-similar Gaussian processes and fields.}
However, verifying these conditions requires in general significant efforts. \\
Given the relation between extremes and the
sojourn time, in this paper we derive also some results for the
mean sojourn time of the $r$th order statistics process
$X_{r:n}$, see below Propositions \ref{P1}, \ref{P2}. Again,
those results will be shown by imposing only conditions on the
generic process $X$. \cL{For more studies on  sojourn time, we refer to  \cite{Berman82, Berman92} and the recent contribution \cite{Pham2013} for stationary Gaussian random fields.}

The rest of this paper is organized as follows: The main
results are displayed in \nesec{sec3} followed then by some  illustrating 
examples. All the proofs are relegated to \nesec{sec5} followed by some concluding remarks.

\COM{The following is a set of conditions suitable for the rest of
the paper:\\ 
$i)$ $\{X(t), t\ge0\}$ is an $\R$-valued self-similar random process with index $\kappa>0$ and $P$-continuous sample
paths;\\ 
$ii)$ $X(1)$ has a 
continuous distribution function (df) with an infinite right endpoint;\\ 
$iii)$ there is a non-increasing, positive function $q=q(u),u>0$  such that
\BQN
\label{eqQ} \lim_{u\to\IF} \frac{1}{ q} =Q \in (0,\IF],
\quad\tilde a = \frac{1}{2\sup_{u \in [0, \IF)}q} >0.
 \EQN
In order to simplify the notation 
we will consider below a
simpler $q$ function, namely we will assume instead throughout
the paper that $i),ii), iii)^*, iv)$ hold with
}

\section{Main Results}\label{sec3} 
In this section we will derive our
main results for the order statistics process $X_{r:n}$ generated by $X$
which was specified in the previous section (see \eqref{def: Orderp}). 
In order to establish the approximations of $\pru$, we need two auxiliary results below on the asymptotic
properties of the mean sojourn time of $r$th order
statistics process $X_{r:n}$ over an increasing threshold $u$, denoted by $L_r(u)=L_r(1; u)$, 
with 
\BQN\label{Def: Lr} 
L_r(s; u) = \int_0^s
\I{ X_{r:n}(t)} \, dt, \quad s\in[0, 1]. 
\EQN 
Throughout in the sequel $G_r(\cdot)$ stands for the df of $X_{r:n}(1)$. 
\BS \label{P1}
If the df $G$ of $X(1)$ is continuous and further $G\in MDA(\Lambda,w)$, then 
\BQNY
\cL{\overline G_r(u) =\frac{n!}{r!(n-r)!}(\overline  G(u))^r (1+o(1)),\quad u\to\IF.}
\EQNY
Moreover, $G_r\in MDA(\Lambda,w_r)$ with $w_r=w_r(u)=
rw(u)$ and 
\BQN \label{Eq: Lr} 
\E{L_r(u)} =\frac1{
\kappa u w_r(u)} \overline G_r(u) (1+o(1)), \quad u\to\IF. 
\EQN
\ES 
More generally, we can obtain below the bounds of the
asymptotic distribution of $L_r(u)$ if conditions A  and B are satisfied by the generic process
$X$. Let $\{\xi_{r:r}(t), t\ge0\}$ be the minimum
process of $r$ independent copies of $\xi$, and define 
\BQN \label{Def: Lambda} 
\Theta_r(x) = \pk{
\int_0^\IF \I{\xi_{r:r}(t)} \, dt >x},\quad x\ge0. 
\EQN 
\BS \label{P2} If condition A
holds for the process $X$, then for each $x\ge 0$
\BQNY 
\liminf_{u\to\IF} \int_x^\IF \frac{\pk{L_r(u)/q>
y}}{\E{L_r(u)/q}}\, dy \ge \Theta_r(x). 
\EQNY 
If additionally condition B holds, then for each
$x\ge0$
\BQNY 
\limsup_{u\to\IF} \int_x^\IF \frac{\pk{L_r(u)/q>
y}}{\E{L_r(u)/q}}\, dy \le \Theta_r(x-). 
\EQNY 
\ES 
{\remark\label{rem0} a) Propositions \ref{P1}, \ref{P2} generalize Proposition 1 and Theorem 1 in \cite{Albin1998} for $r$th order statistics processes, respectively. \\
b) If conditions A, B are satisfied by the generic process $X$, then for any $x>0$ 
\BQNY
\notag \pk{\sup_{t\in[0, 1]}X_{r:n}(t) >u} &\ge &
\max\left(\overline G_r(u), \frac1x\int_0^x\pk{\frac{L_r(u)}{q}>y}\,
dy\right) \\ & \ge & \max\left(\overline G_r(u),
\frac{1-\Theta_r(x-)}x\E{\frac{L_r(u)}{q}}\right). 
\EQNY
}
We will see that the asymptotic inequality above can be reversed under certain tightness conditions C or C$^*$ in \netheo{T3} and \netheo{T4}  below. 
Further by Proposition \ref{P1}, these asymptotic results will depend on 
the limits $\beta_3$ and $\beta_4$ defined as below
\BQN \label{B34} 
\beta_3 = \liminf_{u\to\IF}uqw
\quad \mbox{and}\quad \beta_4 = \limsup_{u\to\IF}uqw. 
\EQN
\BT\label{T3} 
If condition C and $\beta_3=\IF$
are satisfied by the generic process $X$, then 
\BQNY 
 \lim_{u\to\IF} \frac1{\overline G_r(u)}
\pk{\sup_{t\in[0,1]} X_{r:n}(t) >u} =1. 
 \EQNY 
 \ET
\BT\label{T4} 
If  condition A is satisfied by the generic process $X$. If additionally condition C and
$0<\beta_3\le \beta_4<\IF$ hold, or instead conditions  B and C$^*$ hold, then 
\BQNY \lim_{u\to\IF} \frac1{\E{L_r(u)/q}}
\pk{\sup_{t\in[0,1]} X_{r:n}(t) >u} = -\Theta_r^\prime(0)
\EQNY 
where $ -\Theta_r^\prime(0)= \lim_{x\downarrow 0}
(1-\Theta_r(x) )/x$ is positive and finite.
\ET
\COM{{\korr \label{korr1} \netheo{T4} still hold if
condition C and $\beta_3 >0$ are replaced by
conditions B and C$^*$. } }
{\remark \label{rem1} 
a) The approximation of $p_{r,T}(u)= \pk{\sup_{t\in [0,T]} X_{r:n}(t)> u}$ for any $T>0$ is readily obtained by using the self-similarity of 
$X_{r:n}$. \\
b) One may consider the case that the marginal df $G$ belongs to the \FRE  \ and Weibull max-domain attraction, respectively; see, e.g.,
\cite{Albin1990, Albin1998}. \\
c) The suitable choice of $q$ function is crucial for the asymptotic results; see \netheo{T5} below. On the other hand,  one may consider some general function $q(\cdot)$ with minor modifications; see \cite{Albin1998}.}

\section{Examples}\label{sec4}
In this section, we first give two classes of self-similar Gaussian processes, namely bi-fractional Brownian motions and sub-fractional Brownian motions, and then investigate a generalized self-similar skew-Gaussian process to illustrate our main results in \nesec{sec3}. 

{\bf Example 1}: (\emph{Bi-fractional Brownian motions (bi-fBm)}) Let $\{B_{h,k}(t),
t\ge0\}$ with $h\in (0, 1), {k\in(0, 1]}$ be a bi-fBm, i.e.,
a centered self-similar Gaussian process with covariance
function given by 
\BQNY 
\E{B_{h,k}(s)B_{h,k}(t)} =
\frac1{2^k}\Big( \big(s^{2h} + t^{2h}\big)^k -
\abs{s-t}^{2hk}\Big), \quad s, \ t\ge0. 
\EQNY 
In particular, the
bi-fBm $B_{h, 1}$ is the fBm with Hurst index $h$. The Lamperti's 
transformation $\tilde B_{h,k}(t)\equiv e^{-\kappa
t}B_{h,k}(e^t)$ with $\kappa= hk$ is a centered stationary
Gaussian process with covariance function 
\BQNY \E{\tilde
B_{h,k}(0)\tilde B_{h,k}(t)} = 1-2^{-k}t^{2hk} +o(\abs t +
\abs{t}^{2hk}), \quad t\to0. 
\EQNY 
{\bf Example 2}: (\emph{Sub-fractional
Brownian motions (sub-fBm)}) 
Let $\{S_h(t), t\ge0\}$ with $h\in
(0,1)$ be a sub-fBm, i.e.,  a centered self-similar Gaussian process with
covariance given by 
\BQNY 
\E{S_h(s)S_h(t)} = s^{2h}+
t^{2h}-\frac12\Big((s+t)^{2h}+\abs{s-t}^{2h}\Big), \quad s,\, t\ge0. 
\EQNY
The Lamperti's  transformation $\tilde S_{h}(t)\equiv
e^{-\kappa t}S_{h}(e^t)$ with $\kappa= h$ is a centered
stationary Gaussian process with covariance function satisfying
\BQNY \E{\tilde S_h(0)\tilde S_h(t)} = (2-2^{2h-1})\left(1
-\frac{1}{2(2-2^{2h-1})}t^{2h} +o(\abs t +
\abs{t}^{2h})\right),\quad t\to0. 
\EQNY 
Note that the
correlation functions of the
above centered stationary Gaussian processes have regular
varying tails at zero. Using the well-known results for
stationary Gaussian processes (see, e.g., \cite{AlbinC, Pickands1969}) and the Lamperti's Propositions in \cite{Albin1998},
we readily see that \netheo{T4} holds for the two classed of self-similar Gaussian
processes above which will further be extended to be of skew version; see \netheo{T5} below.

{\bf Example 3}: (\emph{Generalized self-similar skew-Gaussian processes}) Let $\{X(t), t\ge 0\}$ be a centered self-similar Gaussian
process with index $\kappa>0$ and covariance function
satisfying
\BQN\label{BoundCov} 
\E{X(1)X(1+t)} = 1+\kappa
t-\mathrm D\abs t^\alpha +o(\abs t+\abs t^\alpha), \quad t\to0 
\EQN 
for some
constants $\alpha\in(0, 2]$ and $\mathrm D>0$.  Denote by $\abs{\boldsymbol{\chi}(t)}
:=\big( \sum_{i=1}^m X_i^2(t) \big)^{1/2}$ with $X_i, 1\le i\le m+1, m\inn$   independent copies of
$X$.  We then define  the generalized skew-Gaussian process
$\zeta$ with $\delta\in [0, 1]$  by
\BQN 
\label{def:
SkewG} \zeta(t) \equiv \delta
\abs{\boldsymbol{\chi}(t)}+\sqrt{1-\delta^2} X_{m+1}(t),
\quad t\ge0. 
\EQN
We note in passing that  condition
\eqref{BoundCov} is satisfied by the bi-fBm and the sub-fBm. In
particular, \eqref{BoundCov} holds with $\mathrm D=1/2, \kappa=\alpha/2$ if
$X(t)=Z(t)$ is a standard fBm with Hurst index $\alpha/2 \in
(0,1]$, i.e., a centered Gaussian process with covariance function satisfying
$$ 
\mathbf{Cov}(Z(s), Z(t))= \frac12 \bigl( s^\alpha+
t^\alpha- \abs{s-t}^\alpha\bigr), \quad s,\, t \ge 0. 
$$
In what follows, utilizing the obtained results in \nesec{sec3}, we investigate the asymptotic results of $\zeta_{r:n}$, the $r$th order
statistics process with the generic process 
$\zeta$. This class of skew-Gaussian processes has received a lot of
attentions from both theoretical and applicable fields; see,
e.g., \cite{AlRawwashS07, AlodatR09, KHorders}.

To formulate our asymptotic results on the process $\zeta$, we need further some notation and condition \eqref{SupboundCov} below. Let $E$ be a unit mean exponential random
variable (rv) which is independent of the standard fBm $Z$ with
Hurst index $\alpha/2\in(0, 1]$. We write further $E_i, i\le r$ and
$Z_i, i\le r$ for independent copies of $E$ and $Z$,
respectively. As for  the $\chi$-process (i.e., $\zeta$ with $\delta=1$) in \cite{Albin1998}, we assume that, there
exists some $h>0$ such that 
\BQN\label{SupboundCov}
\sup_{t\in[\ve, h]}e^{-\kappa t} \E{X(1)X(e^t)} <1,\quad
\ve\in(0,h].
\EQN
\BT\label{T5} Let $\{\zeta(t), t\ge0\}$ be the generalized self-similar skew-Gaussian process given in \eqref{def:
SkewG}. We have\\
a) If \eqref{BoundCov} and \eqref{SupboundCov} hold with
$\alpha\in(0, 1]$, then Propositions \ref{P1}, \ref{P2} and \netheo{T4} hold for the $r$th order statistics process $\zeta_{r:n}$ with 
$w(u) = \max(1,u), q(u)=(\max(1,u))^{-2/\alpha}$ and 
\BQNY \Theta_r(x) =\pk{\int_0^\IF
\I{\min_{1\le i\le r}(\sqrt 2Z_i(\mathrm D^{1/\alpha}t) -\mathrm Dt^\alpha +
E_i-\beta_4\kappa t)}\, dt >x}
\EQNY 
where $\beta_4 = 1$ if
$\alpha =1$, and 0 otherwise.\\ 
b) If \eqref{BoundCov} holds
with $\alpha\in(1, 2]$, then Propositions \ref{P1}, \ref{P2} and \netheo{T4} hold for  $\zeta_{r:n}$  with $w(u) = \max(1, u), q(u)=(\max(1,u))^{-2}$ and $\Theta_r(x) =e^{-\kappa rx}$.
\ET
{\remark 
It is possible to allow that $X_i, 1\le i\le m+1$ in \eqref{def:
SkewG} are not identically distributed, for instance, assuming that the processes $X_i, 1\le i\le m+1$ satisfy condition \eqref{BoundCov} with different $\mathrm{D}_i>0, 1\le i\le m+1$.
}
\COM{ 
It follows easily that $\pk{\zeta_{r:n}(1) >u}=
n!/(r!(n-r)!)\big(\pk{\zeta(1) > u}\big)^r
(1+o(1))$ with (set $0^0:=1$ if $\delta=0, m=1$) \BQN
\pk{\zeta(1) > u} = \delta^{m-1}
\frac{2^{1-m/2}}{\Gamma(m/2)} u^{m-2} \expon{-\frac{u^2}2}
(1+o(1)). \EQN
\BT\label{T6} Let $\{\zeta(t), t\ge0\}$ be
defined by \eqref{def: Orderp} to be satisfied by \ET }
\COM{{\remark\label{rem2} We could relax the assumption of
identical distributions for our results. For instance, suppose
$\zeta^{(j)}_{m, \delta}, j=1, \ldots, n$ have generalized
self-similar skew-Gaussian distributions with index $\kappa$
andcovariance function satisfying \eqref{SkewGcovariance} with
parameters $D_i>0, \alpha_i\in(0, 2]$, then \\ a) \\ b). }
Next,we consider the extremal behavior of the iterated Brownian
motion $\{Z(t), t\ge0\}$, i.e., It follows from Lemma 2.1 in
\cite{ArendarzykD2011} that the distribution function $G$ of
$Z(1)\equaldis X(\abs{Y(1)})$ satisfies that \BQNY 1-G(u) =
\pk{\sqrt{Y(1)}X(1) >u} = C\expon{-3\cdot 2^{-5/3} u^{4/3}}
(1+o(1)), \quad u\to\IF, \EQNY which implies that $G\in
MDA(\Lambda,w)$ with auxiliary function $w = \sqrt[3]{2u}$.}
\section{Proofs}\label{sec5} In this section, we will first
give, by verifying that conditions A and B imposed on $X$ are satisfied by
the order statistics $X_{r:n}$ in turn, the proofs of
Propositions \ref{P1}, \ref{P2}. Then our main results
(Theorems \ref{T3}--\ref{T4}) follow mainly by
the verifications of the tightness conditions C and
C$^*$. We conclude this section with the proof of Theorem
\ref{T5} concerning the generalized 
self-similar skew-Gaussian  processes.

In what follows, we write $\equaldis$ for equality in
distribution while $\todis$ for the convergence in distribution
(or the convergence of finite-dimensional distribution if both
sides of it are random processes). In the sequel, we set 
\BQNY 
c_{n,l}=\frac{n!}{l!(n-l)!}, \quad l=0, \ldots, n\quad{\rm and} \quad u_x=u+\frac xw, \quad x\ge0.
\EQNY
\proofprop{P1} It follows by Lemma 1 in \cite{KHorders} that 
\BQNY
\overline G_r(u) = c_{n,r} (\overline G(u))^r (1+o(1)),\quad u\to\IF
\EQNY
which together with the fact that $G\in MDA(\Lambda, w)$ implies that 
$G_r\in MDA(\Lambda, w_r)$ with $w_r(u)= rw(u)$. Noting further that $\{X_{r:n}(t), t\ge0\}$ is a
self-similar process with index $\kappa$, the claim follows by Proposition 1 in \cite{Albin1998}. \QED

\proofprop{P2} 
It follows from
Lemma 2 in \cite{Albin1998} that $\pk{\xi(t) > x}$ is
continuous at $x=0$ for each $t\in(0, \IF)$. Thus, in view of
Theorem 1 in
\cite{Albin1998}, the lower bound follows if we verify that
\eqref{Eq: condL} holds for the order statistics process
$X_{r:n}$. The main idea of the following proof is from that of Lemma 2 in \cite{KHorders}.

We first consider the proof for $r=n$. By condition A and Lemma 1 in \cite{KHorders}
\BQN\label{eq: condm} 
\lefteqn{\pk{X_{n:n}(1- q t_i) > u, 1\le i\le
m \lvert X_{n:n}(1)> u} }\notag \\ && = \frac{\pk{X_{n:n}(1-
qt_i) > u, 1\le i\le m, X_{n:n}(1)> u}}{\pk{X_{n:n}(1)> u}} \notag\\
&& \to \pk{\xi_{n:n}(t_i) > 0, 1\le i\le m }, \quad
u\to\IF\qquad\qquad\qquad\qquad\qquad\qquad 
\EQN 
establishing the claim for $r=n$. 
The remaining cases that $r<n$ follow
 if we show that
 \BQN\label{CondOrder} 
 \lefteqn{ \pk{X_{r:n}(1- q t_i) > u, 1\le i\le m
\lvert X_{r:n}(1)> u} } \notag \\ && =
\pk{X_{r:r}(1- q t_i) > u, 1\le  i\le
m\lvert X_{r:r}(1)> u} (1+ \Upsilon_r( u))
\EQN 
with $\lim_{u\to\IF}\Upsilon_r( u)=0$ holds uniformly for all $ t_i\in (0, \IF)$ and $m\in \N$.

In the following, we only present the proof for the case that
$r= n-1$ and $m=1$ since the other cases follow by similar
arguments. It follows from Lemma 1 in \cite{KHorders} that 
\BQNY
\overline G_{n-1}(u) = n \pk{X_{(n-1):(n-1)}(1) > u} (1+
o(1)),\quad u\to\IF 
\EQNY 
and 
\BQNY \lefteqn{ \pk{X_{(n-1):n}(1- q t_1) > u,
X_{(n-1):n}(1) > u}} \\ && = \pk{X_{(n-1):n}(1- q t_1) > u\ge
X_{n:n}(1- q t_1), X_{n: n}(1) > u} \\ && \quad+\pk{X_{n:n}(1-
qt_1) > u, X_{(n-1):n}(1) > u\ge X_{n:n}(1)}\\ &&
\quad+\pk{X_{(n-1):n}(1- q t_1) > u\ge X_{n:n}(1- q t_1),
X_{(n-1):n}(1) > u\ge X_{n:n}(1)}\\ && \quad+\pk{ X_{n:n}(1- q
t_1) > u, X_{n: n}(1) > u}\\ && =: I_{1u} + I_{2u} + I_{3u} +
I_{4u}. 
\EQNY 
Since 
$\pk{X_{n}(1- q t_1)\le u, X_{n}(1) > u} \le \overline G(u)=
o(1)
$ 
holds uniformly for all $ t_1\in(0,\IF)$ and large $u$ we have 
\BQNY
I_{1u} &=& n \pk{ \min_{1 \le j \le n-1} X_{j}(1- q t_1)> u,
\min_{1 \le j \le n-1} X_{j}(1)> u } \pk{X_{n}(1- q t_1) \le u,
X_{n}(1) > u}  \\ &=& n \pk{ \min_{1 \le j \le n-1} X_{j}(1-
qt_1)> u, \min_{1 \le j \le n-1} X_{j}(1)> u } o(1), \quad u\to\IF. 
\EQNY
Similarly, $I_{2u} = I_{1u}(1+ o(1)),\ u\to\IF$. \\
Next, we deal with $I_{3u}$ and $I_{4u}$. \cL{Note that 
\eqref{eq: condm} implies that for $k=0, 1, 2 $
\BQNY
\pk{\min_{1\le j \le n-k} X_{j}(1- q t_1) > u, \min_{1\le j'
\le n-k} X_{j'}(1) > u} = (\overline G(u))^{n-k} O(1),
\quad u\to\IF. 
\EQNY
}
Using further the fact that $\pk{X_{n}(1- q t_1)\le u, X_{n}(1) \le u}
= 1+o(1)
$ holds for $ t_1\in(0,\IF)$ and large $u$
\BQNY  I_{3u} &=& \sum_{i, i'=1, \ldots,
n}\pk{\min_{1\le j \le n, j\neq i} X_{j}(1- q t_1) > u,
X_{i}(1-q t_1) \le u, \min_{1\le j' \le n, j'\neq i'} X_{j'}(1) >
u,
X_{i'}(1) \le u}  \\ &=& n \pk{\min_{1\le j \le n-1} X_{j}(1- q
t_1) > u, \min_{1\le j' \le n-1} X_{j'}(1) > u} \pk{X_{n}(1- q
t_1) \le u, X_{n}(1) \le u} \\ && +c_{n, 2}
\pk{\min_{1\le j \le n-2} X_{j}(1- q t_1) > u, \min_{1\le j'
\le n-2} X_{j'}(1) > u} \\ && \times \pk{X_{n-1}(1- q t_1)
\le u, X_{n-1}(1) > u}\pk{X_{n}(1- q t_1) > u, X_{n}(1) \le u} \\ 
&=&n \pk{\min_{1\le j \le n-1} X_{j}(1- q t_1) > u, \min_{1\le j'
\le n-1} X_{j'}(1) > u} (1+ o(1))
\EQNY 
we conclude that $I_{4u} = o(I_{3u}),\ u \to\IF$ and further that 
\eqref{CondOrder} follows for
$r=n-1$ and $m=1$. Thus, we complete the proof of the lower
bound.

Next, we consider the upper bound. By Theorem 1 in
\cite{Albin1998}, we only need to show that condition B is satisfied by the $r$th order process $X_{r:n}$, which can be shown by verifying the following inequality: for sufficiently
large $u$ and some $D>0$ 
\BQN \label{Ineq: cond_r}
\pk{X_{r:n}(1-qt)> u \lvert
X_{r:n}(0) > u} \le D\pk{X(1-qt) > u \lvert X(0)> u}
\EQN
holds locally uniformly with $t\in(0, \IF)$. Clearly, for $r=n$
\BQNY 
\pk{X_{n:n}(1-qt)> u \lvert X_{n:n}(1) > u} =
\Bigl(\pk{X(1-qt) > u \lvert X(1)> u} \Bigr)^n \EQNY and for
$r<n$, similar arguments as for \eqref{CondOrder} yield that
\BQNY 
\pk{X_{r:n}(1-qt) >u\lvert X_{r:n}(1) >u} =
\left(\pk{X(1-qt) > u \lvert X(1) >u}\right)^r
(1+\Upsilon_r(u))
\EQNY 
with some function $\Upsilon_r(u)$ which tends to 0 as $u\to\IF$. 
It thus follows that \eqref{Ineq: cond_r} 
holds uniformly for $t\in(0, \IF)$ and some positive
constant $D$. Therefore, we complete the proof. \QED

\prooftheo{T3} In view of Theorem 4 in \cite{Albin1998}, it
suffices to show the tightness condition C holds also for the order
statistics $X_{r:n}$. To this end, we will show that, for given $\sigma>0$ and small $a\in(0,
\tilde a]$
\BQN\label{vr'} 
\vk{v}_r(a, \sigma)\equiv
\limsup_{u\to\IF} \frac{\pk{\sup_{t\in[0,1]} X_{r:n}(t) >
u_\sigma, \max_{0\le k\le K} X_{r:n}(t_a^u(k)) \le u } } {
\E{L_r(u)/q}+ \overline G_r(u)} \le  D\vk v(a, \sigma)
\EQN
with $\{t_a^u(k)\}$ given by \eqref{t-sequence}, and some positive constant $D$ whose value doesn't depend on $a$ and might change below from line to line. \\  
Letting $q_\sigma= q(u_\sigma)$, we have 
\BQN\label{Bound_ineq}
\lefteqn{ \limsup_{u\to\IF}
\frac{\pk{\sup_{t\in[0,1]} X(t) >u_\sigma}}{\E{L(u)/q}
+\overline G(u)}}\notag \\ &&  \le \limsup_{u\to\IF}\frac{\pk{\sup_{t\in[0,1]}
X(t) > u_\sigma}}{\E{L(u_\sigma)/q_\sigma} +\overline G(u_\sigma)} \limsup_{u\to\IF}\frac{\E{L(u_\sigma)/q_\sigma}
+\overline G(u_\sigma)}{\E{L(u)/q} +\overline G(u)}
\EQN
which, in view of  Proposition 1 and Theorem 3 in
\cite{Albin1998}, is bounded. Using
further Proposition \ref{P1}, $\beta_3>0$ and the inequality $(x+y)^n\le 2^n(x^n +y^n), \ x, y>0$
\BQNY
\lefteqn{ \pk{\sup_{t\in[0,1]} X_{n:n}(t) > u_\sigma,
\max_{0\le k\le K} X_{n:n}(t_a^u(k)) \le u } } \\ && \le \pk{
\sup_{t\in[0,1]} X_j(t) >u_\sigma, 1\le j\le n,
\cup_{i=1}^n\{\max_{0\le k\le K} X_i(t_a^u(k)) \le u \}} \\
&& \le \sum_{i=1}^n \pk{ \sup_{t\in[0,1]} X_i(t) >u_\sigma,
\max_{0\le k\le K} X_i(t_a^u(k)) \le u }
\Biggl(\pk{\sup_{t\in[0,1]} X(t) >u_\sigma}\Biggr)^{n-1} \\
&& \le D\left(\E{\frac{L(u)}{q}} + \overline G(u)\right)^n \vk{v}(a, \sigma) \\
&&\le D\left(\E{\frac{L_n(u)}{q}} + \overline G_n(u)\right)\vk{v}(a, \sigma), \quad u\to\IF
\EQNY
establishing the proof of \eqref{vr'} for $r=n$.\\ 
Next, we present only the proof for $r=n-1$ since the other cases follow by similar
arguments. \\ 
Note that $\vk v(a, \sigma)$ is non-negative and finite for small $a>0$ and $\sigma>0$. We have 
 \BQNY \lefteqn{
\pk{\sup_{t\in[0,1]} X_{(n-1):n}(t) > u_\sigma, \max_{0\le k\le
K} X_{(n-1):n}(t_a^u(k)) \le u } } \\ && \le n\pk{
\sup_{t\in[0,1]} X_j(t) >u_\sigma, 1\le j\le n-1, \cup_{i, j=1,
\ldots,n}\{\max_{0\le k\le K} X_i(t_a^u(k)) \le u, \max_{0\le
k\le K} X_j(t_a^u(k)) \le u \}} \notag\\ && \le n\sum_{i,
j=1,\ldots,n-1} \pk{ \sup_{t\in[0,1]} X_i(t) >u_\sigma,
\max_{0\le k\le K} X_i(t_a^u(k)) \le u } \notag\\ &&
\quad\times \pk{ \sup_{t\in[0,1]} X_j(t) >u_\sigma, \max_{0\le
k\le K} X_j(t_a^u(k)) \le u } \Biggl(\pk{\sup_{t\in[0,1]} X(t)
>u_\sigma}\Biggr)^{n-3} \notag \\ && \quad+
2n\sum_{i=1,\ldots,n-1, j=n} \pk{ \sup_{t\in[0,1]} X_i(t)
>u_\sigma, \max_{0\le k\le K} X_i(t_a^u(k)) \le u }
\Biggl(\pk{\sup_{t\in[0,1]} X(t) >u_\sigma}\Biggr)^{n-2}
\notag\\ && \le D n\big( c_{n-1,2}\vk{v}(a,
\sigma)+2(n-1)\big) \Biggl(\E{\frac{L(u)}{q}}
+\overline G(u)\Biggr)^{n-1} \vk{v}(a, \sigma)\notag\\ && \le
D\left(\E{\frac{L_{n-1}(u)}{q}} + \overline G_{n-1}(u)\right)\vk{v}(a,\sigma)
\EQNY
where the third inequality follows by \eqref{Bound_ineq} and condition C. 
It follows thus that \eqref{vr'}
holds for $r=n-1$. We complete the proof. \QED

\prooftheo{T4} We first present the proof under conditions A, C and $0<\beta_3\le \beta_4<\IF$.
To this end, we will show conditions A, B and C are satisfied by the $r$th order statistics process $X_{r:n}$ in turn.  \\
By Proposition \ref{P1},  $G_r\in MDA(\Lambda,w_r)$. 
Further, it follows from the proved \eqref{eq: condm} and \eqref{CondOrder} that conditin A is satisfied by the $r$th process $X_{r:n}$ and the limit process $\xi_{r:r}$. In view of Proposition 2 in \cite{Albin1998}, 
$\beta_3 >0$ imply that condition B holds for the order process $X_{r:n}$.  Finally, the proved \eqref{vr'} shows that condition  C holds also for the
order process $X_{r:n}$. \\
Therefore, in view of Theorems 5, 6 in
\cite{Albin1998}, we have 
\BQNY 
\lefteqn{ \limsup_{u\to\IF}
\frac1{\E{L_r(u)/q}} \pk{\sup_{t\in[0,1]} X_{r:n}(t) >u} }\\ && \le
\liminf_{x\downarrow 0} \frac{1-\Theta_r(x)}x \le
\limsup_{x\downarrow 0} \frac{1-\Theta_r(x)}x \\
&&\le
\liminf_{u\to\IF} \frac1{\E{L_r(u)/q}} \pk{\sup_{t\in[0,1]}
X_{r:n}(t) >u}
\EQNY 
with $\Theta_r(\cdot)$ given by \eqref{Def: Lambda}. \\
Noting that conditions A, B are satisfied by $X_{r:n}$ , we have by Theorem 2 in \cite{Albin1998} 
that the limit $-\Theta_r'(0)$ is positive and finite. Hence the claim follows. 

 Next, we consider instead condition C$^*$ holds also  for the order statistics
process $X_{r:n}$, i.e., with  the involved constants given as
in condition C$^*$
\BQN\label{condC*O} 
\pk{X_{r:n}(1-qt)
>u_{\lambda+v}, X_{r:n}(1) \le u_v} \le D^*
t^d\lambda^{-b}\overline G_r(u) 
\EQN
for all $u$ large and all $0< t^\rho\le \lambda \le \lambda_0, v\ge0, D^*>0$.  

It follows by the self-similarity of  $X$ that $\pk{X(1-qt)
>u_{\lambda+v}}\le \overline G(u)
$ holds uniformly for
$1-qt\in(0, 1)$. We have
\BQNY 
\lefteqn{ \pk{X_{n:n}(1-qt)
>u_{\lambda+v}, X_{n:n}(1) \le u_v}}\\ &&
=\pk{X_i(1-qt) >u_{\lambda+v}, 1\le i\le n,
\cup_{j=1,\ldots, n} \left\{X_j(1) \le u_v\right\}} \\ && \le
n\pk{X(1-qt) >u_{\lambda+v}, X(1) \le u_v}
\left(\pk{X(1-qt) >u_{\lambda+v}}\right)^{n-1}\\ && \le
n Dt^d \lambda^{-b} \left( \overline G(u)\right)^n = n Dt^d \lambda^{-b}\overline G_n(u) 
\EQNY
establishing the proof of \eqref{condC*O} for $r=n$. \\
Next, we only deal with the case $r=n-1$ since the other cases
follow by similar arguments. 
\BQNY 
\lefteqn{
\pk{X_{(n-1):n}(1-qt) >u_{\lambda+v}, X_{(n-1):n}(1)
\le u_v}}\\ && = n\pk{X_i(1-qt) >u_{\lambda+v},
1\le i\le n-1, \cup_{i,j=1,\ldots, n} \left\{X_i(1) \le u_v, X_j(1)
\le u_v\right\}} \\ && \le n\sum_{i, j=1,\ldots,
n-1}\left(\pk{X(1-qt) >u_{\lambda+v}, X(1) \le u_v}\right)^2 \left(\pk{X(1-qt)
>u_{\lambda+v}}\right)^{n-3}\\ && \quad + 2
n\sum_{i=1,\ldots, n-1, j=n}\pk{X_i(1-qt)
>u_{\lambda+v}, X_i(1) \le u_v}\left(\pk{X(1-qt) >u_{\lambda+v}}\right)^{n-2} \\
&& \le nD\big(c_{n-1,2}+2(n-1)\big)t^d \lambda^{-b}
\left(\overline G(u)\right)^{n-1} \\
&&\le D^*t^d
\lambda^{-b} \overline G_{n-1}(u),\quad u\to\IF
\EQNY holds for all $0 < t^\rho\le \lambda \le \lambda_0, v\ge0$ and $D^*$ a
positive constant. Thus \eqref{condC*O}  holds for $r=n-1$. Consequently,  the desired result follows from Corollary 1 in \cite{Albin1998}. \QED

\prooftheo{T5} Without loss of generality, we assume that $\mathrm{D}=1$. In the following, we denote 
the Lamperti's  associated stationary process $\tilde \zeta(t) = e^{-\kappa t}\zeta(e^t), \, t>0$. We will show below that conditions A, B and C are satisfied by the generic process $\zeta$ in turn, and thus Propositions \ref{P1}, \ref{P2} and \netheo{T4} follow. \\
a)  We start below with the verification of condition A. It follows from Lemma 5
in \cite{KHorders} that the marginal distribution $G$ of $\zeta(1)$ satisfies that $G\in
MDA(\Lambda,w)$ with auxiliary function $w(u)=\max(1,u)$. Letting $q=q(u)=(\max(1, u))^{-2/\alpha}, \alpha\in(0,1]$, we have $\beta_3=\beta_4 =1$ if
$\alpha=1$, 0 otherwise. Further, it follows from Lemma 6 in \cite{KHorders} that the f.d.d. of  $\left\{w(u)(\tilde \zeta(-qt)
-u)\Big\lvert (\tilde \zeta(0) >u),\, t>0\right\}$ converges to
those of $\sqrt 2Z(t) -t^\alpha + E$. Thus, in view of
Proposition 9 (ii) in \cite{Albin1998}, condition A holds with
$\xi(t) \stackrel{d}{=}\sqrt 2Z(t) -t^\alpha + E -\beta_4
\kappa t.$ \\
Next, it follows from Lemma 7 in
\cite{KHorders} that there exists some positive constant $K_p$ and $p > m$ 
such that 
\BQNY 
\pk{\tilde \zeta(qt) > u  \lvert
\tilde \zeta(0) > u} \le  \left\{
\begin{array}{ll}K_p t^{-\alpha p/2}, & qt\in (0, \epsilon], \\
K_p u^{m-1-p}, &
qt\in (\epsilon, T] 
\end{array} 
\right. 
\EQNY 
which together
with Proposition 7 in \cite{Albin1998} implies that condition B
holds. \\
Finally, by Lemma 8 in \cite{KHorders}, there exist some
positve constants $D^*, p, \lambda_0$ and $ d>1$ such that
\BQN\label{L3.7} \pk{\tilde \zeta(qt) > u_\lambda, \tilde \zeta(0) \le u} \le
D^*t^d\lambda^{-p} \pk{\tilde \zeta(0)>u} \EQN for
$0< t^{\alpha/2} \le \lambda \le \lambda_0$ and large $u$. Thus, it
follows from Proposition 2 in \cite{Albin1992a} and Proposition
8 in \cite{Albin1998} that condition C is satisfied
by the generalized self-similar skew-Gaussian process $\zeta$. Consequently, we complete the proof of a).\\ 
b) Letting $w=w(u)$ as in a) and $q=q(u)= (\max(1,u))^{-2}$, we have $\beta_4=1$. It follows from Lemma 6 in \cite{KHorders} that the f.d.d. of $\left\{w(u)(\tilde \zeta(-qt) -u)\Big \lvert (\tilde
\zeta(0) >u), \, t>0\right\}$ converges to those of $E$. Thus, in
view of Proposition 9 (ii) in \cite{Albin1998}, condition A
holds with $\xi(t) \stackrel{d}{=}E - \kappa t.$

Next, we verify that conditions B and C are satisfied by the process $\zeta$.  Noting that $\beta_3 =1>0$, condition B follows from
Proposition 2 in \cite{Albin1998}.
Using \eqref{BoundCov} we have 
$\E{\tilde X(0)\tilde X(t)} \ge 1-2\abs t, \, t\to 0$ with $\tilde X(t)= e^{-\kappa t}X(e^t)$.
It follows further from the arguments of Lemma 8 in \cite{KHorders} that \eqref{L3.7} holds. Therefore, in view of  Proposition 2 in \cite{Albin1992a} and Proposition
8 in \cite{Albin1998} condition C hold. 

Consequently, the desired result follows with an elementary calculation to give that $\Theta_r(x) = e^{-\kappa rx}, \, x >0$.
\QED

 \section{Conclusions}\label{sec6}
In this paper, we study  approximations of the supremum of order statistics of self-similar processes. The methodology developed in the seminal paper 
\cite{Albin1998} is powerful for the establishment of our main results.  As mentioned
therein, results on extremes for a self-similar process do not on their own imply results for Lamperti's  associated stationary process or vice versa.  Albin's theory on
extremes of self-similar processes is based on several conditions
(referred to as Albin's conditions) which precisely determine the
asymptotic behaviour of supremum of self-similar processes. As shown by 
the illustrating example of the generalized self-similar skew-Gaussian process, checking the conditions therein requires generally a lot of technical details. Results
for the supremum of self-similar or more general Gaussian processes
can be also derived by utilizing the double-sum method (see, e.g.,
\cite{Pit96}). The importance of Albin's theory lies on the fact
that no Gaussianity assumption is needed rendering very general
results for general self-similar processes.  These theoretical results would be very useful in many applied-oriented fields. 

\textbf{Acknowledgments.}
Partial support from the Swiss National Science Foundation Project 200021-140633/1  and the project RARE -318984
 (an FP7 Marie Curie IRSES Fellowship) are kindly acknowledged. 

\def\polhk#1{\setbox0=\hbox{#1}{\ooalign{\hidewidth
  \lower1.5ex\hbox{`}\hidewidth\crcr\unhbox0}}}

\end{document}